\documentclass{llncs}

\usepackage{amsmath,amsfonts,amssymb,amstext,empheq}
\usepackage{booktabs}
\usepackage{graphicx}
\usepackage{hyperref}

\newcommand{\Cbar}{\bar C}

\newcommand{\Atilde}{\tilde A}
\newcommand{\atilde}{\tilde a}

\newcommand{\Ctilde}{\tilde C}
\newcommand{\Etilde}{\tilde E}
\newcommand{\etilde}{\tilde e}
\newcommand{\abar}{\bar a}

\newcommand{\Aa}{{\mathcal A}}

\newcommand{\Xx}{{\mathcal X}}
\newcommand{\Tt}{{\mathcal T}}

\newcommand{\diag}{\mathrm{diag}}
\newcommand{\Diag}{\mathrm{Diag}}
\newcommand{\inprod}[2]{{\langle #1,#2 \rangle}} 
\newcommand{\st}{\mathrm{s.t.}~~}
\newcommand{\twovector}[2]{\begin{pmatrix}#1\\#2\end{pmatrix}}
\newcommand{\trace}{\mathrm{trace}}

\spnewtheorem{obs}[theorem]{Observation}{\bfseries}{\itshape}

\title{Exact Solution Methods for the $k$-item Quadratic Knapsack Problem}
\titlerunning{An Exact Solution Method for $k$-QKP}
\author{Lucas L\'etocart\inst{1} and Angelika Wiegele\inst{2}}
\authorrunning{L. L\'etocart and A. Wiegele}
\institute{Universit\'e Paris 13, Sorbonne Paris Cit\'e, LIPN, CNRS, (UMR 7030), 93430~Villetaneuse, France, \email{lucas.letocart@lipn.univ-paris13.fr}
\and Alpen-Adria-Universit\"at Klagenfurt, 9020~Klagenfurt am W\"orthersee, Austria, \email{angelika.wiegele@aau.at}}

\begin{document}

\maketitle
\begin{abstract}
The purpose of this paper is to solve the 0-1 $k$-item quadratic knapsack
problem $(kQKP)$, a problem of maximizing a quadratic function subject to two 
linear constraints. We propose an exact method based on semidefinite optimization.
The semidefinite relaxation used in our approach includes simple rank one constraints, which can be handled efficiently by interior point methods. 
Furthermore, we strengthen the relaxation by polyhedral constraints and obtain approximate solutions to this semidefinite problem by applying a bundle method. 
We review other exact solution methods and compare all these approaches by experimenting with instances of various sizes and densities.

\end{abstract}

\section{Introduction}
The 0-1 $k$-item quadratic knapsack problem consists of maximizing a quadratic objective function subject to a linear capacity 
constraint with an additional equality cardinality constraint:

			$$(kQKP) \left\{
			\begin{array}{lllll}
	 		\max &f(x)= \sum_{i=1}^{n}\sum_{j=1}^{n}c_{ij}x_ix_j \\
				\mathrm{s.t.}&
				\sum_{j=1}^n a_j {x}_j \leq b\ \ \ \ &(1)\\
			&\sum_{j=1}^n x_j = k\ \ \ \ \ \ \ &(2)\\
				&x_{j}\in \{0,1\} \ \ \ \  j = 1, \ldots, n
			\end{array}
			\right.$$

\noindent where $n$ denotes the number of items, and all the data, $k$ (number of items to be filled in the knapsack), 
$a_j$ (weight of item $j$), $c_{ij}$ (profit associated with the selection of items $i$ and $j$) and $b$ (capacity of the 
knapsack) are nonnegative integers. Without loss of generality, matrix $C=(c_{ij})$ is assumed to be symmetric.\\
Moreover, we assume that $\max_{j=1,\ldots,n} a_j \leq b < \sum_{j=1}^n a_j$ in order to avoid either trivial solutions or
 variable fixing via constraint (1). Let us denote by $k_{max}$ the largest number of items which could be filled in the
 knapsack, that is the largest number of the smallest $a_j$ whose sum does not exceed $b$. We can assume that 
$k \in \{2,\ldots,k_{max} \}$, where $k_{max}$ can be found in $O(n)$ time \cite{BAZ80,FAY82}. Otherwise, either the value of 
the problem is equal to $\max_{i=1,\ldots,n} c_{ii}$ (for $k=1$), or the domain of $(kQKP)$ is empty (for $k > k_{max}$).\\

$(kQKP)$ is an NP-hard problem as it includes two classical NP-hard subproblems, the $k$-cluster 
problem \cite{BIL05} by dropping  constraint (1), and the quadratic knapsack
problem \cite{PIS07}  by dropping  constraint (2).
Even more, the work of Bhaskara et. al~\cite{BhChGu:12} indicates that
approximating $k$-cluster within a polynomial factor might be a harder problem
than Unique Games. 
Rader and Woeginger~\cite{RW02} state negative results concerning the approximability of QKP if
negative cost coefficients are present.

Applications of $(kQKP)$ cover those found in previous references for $k$-cluster or classical quadratic knapsack problems (e.g., task assignment problems in a
 client-server architecture with limited memory), but also multivariate linear regression and portfolio selection. Specific 
heuristic and exact methods including branch-and-bound and branch-and-cut with surrogate relaxations have been designed for 
these applications (see, e.g., \cite{BS09,BIE96,BOL09,MES07,SLK08}).\\


The purpose of this paper is twofold. 
\begin{enumerate}
\item We introduce a new algorithm for solving $(kQKP)$ and 
\item we briefly review other state of the art methods and compare the methods
  by running numerical experiments.
\end{enumerate}
Our new algorithm consists of a branch-and-bound framework using
\begin{itemize}
\item a combination of a semidefinite relaxation and polyhedral cutting planes to obtain tight
  upper bounds and
\item fast hybrid heuristics~\cite{LET14} for computing high quality lower
  bounds.
\end{itemize}

This paper is structured as follows. In Sect.~\ref{sec:sdp} a semidefinite
relaxation is derived, followed by a discussion of solving the
semidefinite problems in Sect.~\ref{sec:solve}. The relaxation is used inside
a branch-and-bound framework, the various components of this branch-and-bound
algorithm are discussed in Sect.~\ref{sec:bab}. Other methods for solving
$(kQKP)$ and numerical results are presented in Sect.~\ref{sec:numerical}, and Sect.~\ref{sec:conclusion}
concludes.

\paragraph{Notation.} We denote by $e$ the vector of all ones of appropriate size. $\diag(X)$ refers to diagonal of $X$ as a vector and $\Diag(v)$ is the diagonal matrix having diagonal $v$.

\section{A semidefinite relaxation of $(kQKP)$}\label{sec:sdp}
In order to develop a branch-and-bound algorithm for solving $(kQKP)$ to optimality we aim in finding strong upper bounds. 
Semidefinite optimization proved to provide such strong bounds, see e.g.~\cite{ReSo:07,ReRiWi:10,AnGhHuLiWi:13}.

A straightforward way to obtain a semidefinite relaxation is the
following. Express all functions involved as quadratic functions,
i.e. functions in $xx^{t}$, replace the product $xx^{t}$ by a matrix $X$ and
get rid of non-convexities by relaxing $X=xx^{t}$ to $X\succeq xx^{t}$.

Hence, we apply the following changes:
\begin{itemize}
\item Replace the constraint $e^{t}x = k$ by the constraint $(e^{t}x-k)^{2}=0$. 
\item As for the capacity constraint, define $b'$ to be the sum
of the weights of the $k$ smallest items. Clearly, $b'\leq a^{t}x$ is a valid
constraint for $(kQKP)$. Combining
this redundant constraint with the capacity constraint we obtain 
$(b'-a^{t}x)(b-a^{t}x) \leq 0$.
\item Transform the problem to a $\pm 1$ problem by setting $y=2x-e$.
\item Relax the problem by relaxing $Y=yy^{t}$ to $Y\succeq yy^{t}$, i.e., dropping the constraint $Y$ being of rank one. 
\end{itemize}
This procedure yields the following semidefinite problem:
\begin{equation}\tag{${SDP_{1}}$}\label{sdp1}
\begin{split}
\max~~ & \inprod{\Ctilde}{Y} \\ 
\st  & \diag(Y)=e\\
     & \inprod{\Etilde}{Y} = 0\\
     & \inprod{\Atilde}{Y} \le (b-b')^{2} \\
     & Y \succeq 0 
\end{split}
\end{equation}
with $\Etilde = \etilde\etilde^{t}$, $\etilde = \twovector{n-2k}{e}$, 
$\Atilde = \atilde\atilde^{t}$, $\atilde = \twovector{a^{t}e-(b+b')}{a}$, and appropriate $\Ctilde$.

\begin{obs}
\eqref{sdp1} has no strictly feasible point and thus Slater's condition does not hold.
\end{obs}
\begin{proof}
Note that $\inprod{\Etilde}{Y}=\etilde^{t}Y\etilde=0$ together with $Y\succeq 0$ implies
$Y$ being singular and thus every feasible solution is singular. 
\qed \end{proof}
Observe that $\etilde = \twovector{n-2k}{e}$ is an eigenvector to the eigenvalue 0 of every feasible $Y$.
Now consider matrix $V=\twovector{\frac{1}{2k-n}e^{t}}{I_{n}}$. 
$V$ spans the orthogonal complement of the span of eigenvector $\etilde$.
Set $Y=VXV^{t}$ to ``project out'' the
0-eigenvalue and consider the $n\times n$ matrix $X$ instead of the $n+1 \times n+1$ matrix $Y$. 
The relationship between $X$ and $Y$ is simply given by
\begin{equation*}
Y = VXV^{t} 
= \left(\begin{array}{cc} 
\frac{1}{(2k-n)^{2}}e^{t}Xe & \frac{1}{2k-n}(Xe)^{t}\\
 \frac{1}{2k-n}Xe & X\\
\end{array}\right).
\end{equation*}
Looking at the effect of the constraints of~\eqref{sdp1} on matrix $X$, we derive the following conditions.
\begin{itemize}
\item From $\diag(Y)=e$ we obtain the constraints 
  \begin{align*}
    e^{t}Xe & = (2k-n)^{2}\\
    \diag(X) &= e
  \end{align*}

\item The left-hand side of constraint $\inprod{\Etilde}{Y}=0$ translates into
  \begin{equation*}
    \inprod{\Etilde}{Y} = \inprod{\Etilde}{VXV^{t}} = \inprod{V^{t}\Etilde V}{X} = \inprod{0}{X} = 0
  \end{equation*}
  and the constraint becomes obsolete.

\item Constraint $\inprod{\Atilde}{Y} \le (b-b')^{2}$ yields the following.
  \begin{align*}
    \inprod{\Atilde}{Y} &= \inprod{\Atilde}{VXV^{t}} = \inprod{V^{t}\Atilde V}{X}=\\
    &=\inprod{(\frac{a^{t}e-(b+b')}{2k-n}e+a)(\frac{a^{t}e-(b+b')}{2k-n}e+a)^{t}}{X}
  \end{align*}
  Hence,
  \begin{align*}
    \inprod{(\frac{a^{t}e-(b+b')}{2k-n}e+a)(\frac{a^{t}e-(b+b')}{2k-n}e+a)^{t}}{X}
    \le (b-b')^{2} 
  \end{align*}
\end{itemize}
Defining $\abar = (\frac{a^{t}e-(b+b')}{2k-n}e+a)$ we finally obtain

\begin{equation}\tag{$SDP$}\label{sdp}
\begin{split}
\max~~ & \inprod{\Cbar}{X} \\ 
\st  & \diag(X)=e\\
     & \inprod{E}{X} = (2k-n)^{2}\\
     & \inprod{A}{X} \le (b-b')^{2} \\
     & X \succeq 0
\end{split}
\end{equation}
where $E=ee^{t}$, $A=\abar\abar^{t}$, and appropriate cost matrix $\Cbar$. 

\paragraph{Strengthening the relaxation.}
Since we derived a relaxation from a problem in $\pm 1$ variables, we can
further tighten the bound by adding the well known {\em triangle
  inequalities} to the semidefinite relaxation~\eqref{sdp}. These are
for any triple $1 \le i < j < k \le n$:
\begin{equation}\label{tri}
\begin{split}
 x_{ij}+x_{ik}+x_{jk} &\ge -1\\
-x_{ij}-x_{ik}+x_{jk} &\ge -1\\
-x_{ij}+x_{ik}-x_{jk} &\ge -1\\
 x_{ij}-x_{ik}-x_{jk} &\ge -1
\end{split}
\end{equation}
For several problems formulated in $\pm 1$ variables adding these constraints
significantly improves the bound, see e.g.~\cite{ReRiWi:10}. The set of
matrices satisfying all triangle-inequalities is called the {\em metric
  polytope} and is denoted by $MET$. Thus, the strengthend semidefinite relaxation
reads
\begin{equation}\tag{$SDP_{MET}$}\label{sdp-met}
\begin{split}
\max~~ & \inprod{\Cbar}{X} \\ 
\st  & \diag(X)=e\\
     & \inprod{E}{X} = (2k-n)^{2}\\
     & \inprod{A}{X} \le (b-b')^{2} \\
     & X \in MET\\
     & X \succeq 0
\end{split}
\end{equation}

\section{Solving the Semidefinite Relaxations}\label{sec:solve}

\subsection{Solving the Basic Relaxation \eqref{sdp}}
The most prominent methods for solving semidefinite optimization problems are
interior point methods. The interior point method is an iterative algorithm where
in each iteration Newton's method is applied in order to compute new search directions.

Consider the constraints $\Aa(X)=(\vdots)$ with $\Aa(X)=
\left(\begin{array}{c}
\langle A_{1},X \rangle\\ 
\langle A_{2},X \rangle\\
\vdots \\
\langle A_{m},X \rangle
\end{array} \right)$. 
In each iteration we determine a search direction $\Delta y$ ($y$ are variables in the dual semidefinte problem) by solving the system $M \Delta y = rhs$ where
\[m_{ij} = \trace(Z^{-1}A_{j} XA_{i}).\]
$Z$ denotes the (positive definite) matrix variable of the dual semidefinite program.

Forming this system matrix requires $O(mn^{3}+m^{2}n^{2})$ steps and is among the most
time-consuming operations inside the interior point algorithm. (The other time-consuming steps are
maintaining positive definiteness of the matrices $X$ and $Z$ and linear algebra
operations such as forming inverse matrices.)

The primal-dual pair of \eqref{sdp} in variables $(X,s,y,Z,t)$ is given as follows.
\begin{align*}
\max \big\{ &\inprod{\Cbar}{X} \\ 
 & \st~ \diag(X)=e, \inprod{E}{X} = (2k-n)^{2}, \inprod{A}{X} + s = (b-b')^{2}, X \succeq 0, s\ge 0\big\}
\end{align*}
\begin{align*}
\min \big\{ &e^{t}y_{1:n} + (n-2k)^{2}y_{n+1} + (b-b')^{2}y_{n+2} \quad\\
& \st~ \Diag(y_{1:n}) + y_{n+1}E + y_{n+2} A - Z = \Cbar, y_{n+2} - t = 0, \ Z\succeq 0, t\ge 0\big\}
\end{align*}
Hence, the set of constraints is rather simple and 
the system matrix $M$ reads
\[
\left(\begin{array}{ccc}
 Z^{-1}\circ X           & \diag(Z^{-1} E    X)          & \diag(Z^{-1}{A}X)\\                          
\diag(Z^{-1} E    X)^{t}& \inprod{ E   }{Z^{-1} E    X} & \inprod{ E   }{Z^{-1}{A}X}\\
\diag (Z^{-1}A X)^{t} & \inprod{A}{Z^{-1} E  X} & \inprod{A}{Z^{-1}A X} + \frac{s}{t}
\end{array}\right).
\]
Even more, all data matrices have rank~one which can be exploited when computing
the inner products, e.g.,
\[
\inprod{A}{Z^{-1}A X} = \trace (\abar\abar^{t} Z^{-1}\abar\abar^{t} X) 
= (\abar^{t} Z^{-1}\abar)(\abar^{t} X\abar)
\]

Thus, the computation of the inner products of the matrices simplifies and computing the system matrix can be reduced from $O(mn^{3}+m^{2}n^{2})$ to $O(mn^{2}+m^{2}n)$. And since $m=n+2$ in our case, we end up with $O(n^{3})$.

Hence, \eqref{sdp} can be solved efficiently by interior point methods.

\subsection{Solving the Strengthened Relaxation \eqref{sdp-met}}
Problem~\eqref{sdp-met} has a considerably larger number of
constraints than~\eqref{sdp}. Remember that $X\in MET$ is described by $4{n \choose 3}$ linear
inequalities and thus solving~\eqref{sdp-met}  by interior point
methods is intractable. An alternative has been proposed in~\cite{FiGrReSo:06}. Therein the
concept of {\em bundle methods} is used, in order to obtain an approximate optimizer
on the dual functional and thus getting a valid upper bound on~\eqref{sdp-met},
leading to a valid upper bound on~$(kQKP)$.

Bundle methods have been developed to minimize nonsmooth convex functions. To
characterize the problem to be solved, an oracle has to be supplied that evaluates the function 
at a given point and computes an $\epsilon$-subgradient. The set of points,
function values, and subgradients is collected in a ``bundle'', which is used
to construct a cutting plane model minorizing the function to be minimized. 
By doing a sequence of {\em descent steps} the cutting plane model is refined and one gets closer to the minimizer of the function. 

We will apply the bundle method to minimize the dual functional of~\eqref{sdp-met}. Let
\[\Xx = \{ X\succeq 0 \colon \diag(X)=e,\ \inprod{E}{X} = (2k-n)^{2},\ \inprod{A}{X}
\le (b-b')^{2}\}\]
i.e., the feasible region of~\eqref{sdp}. We introduce the dual functional
\begin{equation}\label{sdp-met-dual}
\begin{split}
f(\gamma) &= \max_{X\in \Xx} \{ \inprod{\Cbar}{X} + \gamma^{t}(e-\Tt(X) \}\\
 &= e^{t}\gamma + \max_{X\in \Xx}  \inprod{\Cbar - \Tt^{t}(\gamma)}{X} 
\end{split}
\end{equation}
where $\Tt(X) \le e$ denotes the triangle inequalities~\eqref{tri}.
Minimizing $f(\gamma)$ over $\gamma \ge 0$ gives a
valid upper bound on~$(kQKP)$. In fact, any ${\tilde \gamma} \ge 0$ gives a valid upper bound
\[ z^{*} = \min_{\gamma \ge 0} f(\gamma) \le f({\tilde \gamma}) \mbox{ for any } {\tilde \gamma} \ge 0.\]
Since we use this bound inside a branch-and-bound framework, this allows us to stop early and prune a node as soon as $f(\tilde \gamma)$ is smaller than some known lower bound.
Furthermore, we do not rely on getting to the optimum. We will stop once we are ``close'' to optimum and branch, rather than investing time in dropping the bound by a tiny number.

Evaluating function~\eqref{sdp-met-dual} (the most time consuming step in the bundle method) amounts in solving~\eqref{sdp} (with
varying cost matrix), which can be done efficiently as discussed in the previous section. 
Having the maximizer $X^{*}$ of~\eqref{sdp}, i.e. the function evaluation, a subgradient is given by $g^{*}=e-\Tt(X^{*})$.

\paragraph{Dynamic version of the bundle method.} The number of variables $\gamma$ in \eqref{sdp-met-dual} is $4{n \choose 3}$. 
This number is substantially larger than the dimension of the problem and we are interested only in those inequalities that are likely to be active at the optimum.
Thus, we do not consider {\em all} triangle inequalities but work with a subset that is updated on a regular basis, say every fifth descent step. 
The update consists of
\begin{enumerate}
\item adding the $m$ inequalities being most violated by the current iterate $X$ and
\item removing constraints with $\gamma$ close to $0$ (an indicator for an inactive constraint).
\end{enumerate}
In this way we are able to efficiently run the bundle algorithm by keeping the size of the variable vector $\gamma$ reasonably small.

\section{Branch and Bound}\label{sec:bab}
We develop an exact solution method for solving $(kQKP)$ by designing a
branch-and-bound framework using relaxation~\eqref{sdp-met} discussed above for
getting upper bounds.

The remaining tools of our branch-and-bound algorithm are described in this
section.

\subsection{Heuristics for obtaining Lower Bounds}\label{sec:heuristics}
We use two heuristics to obtain a global lower bound inside our algorithm: one that is executed at the root node and another one that is called at each other node in the branch-and-bound tree.

As a heuristic method at the root node we chose the primal heuristic denoted by $H_{pri}$ in~\cite{LET14}, which is an adaption of a well-known heuristic 
developed by Billionnet and Calmels~\cite{BC96} for the classical quadratic knapsack problem (QKP). This primal heuristic combines a greedy algorithm with local search.

At each node of the branch-and-bound tree, we apply a variable fixation heuristic inspired from $H_{sdp}$~\cite{LET14}. This heuristic method uses the solution of the semidefinite relaxation
obtained at each node, it fixes variables under some treshold $\epsilon >0$ to zero and applies the primal heuristic over the reduced problem. It updates the solution by performing a fill-up and exchange procedure 
over the unreduced problem. This procedure iterates, increasing $\epsilon$ at each iteration, until the reduced problem is empty.

Both heuristics, the primal and the variable fixation one, are very fast and take only hundredths of a second for sizes of our interest.

\subsection{Branching Rule and Search Strategy}
As a branching variable we choose the ``most fractional'' variable, i.e.,
$v=\mathrm{argmin}_{i}|\frac{1}{2} - x_{i}|$. The vector $x$ is extracted from matrix $X$ given by the semidefinite relaxation.

We traverse the search tree in a {\em best first search} manner, i.e., we
always consider the node in the tree having the smallest upper bound.


\subsection{Speed up for small $k$}
Whenever $k$, the number of items to be filled in the knapsack, is small, a branch-and-prune algorithm is triggered in order to speed-up the approach. No relaxation is performed at each node of the branch-and-prune tree and 
a fast depth first search strategy, in priority fixing variables to one, is implemented. We only check the feasibility of the current solution through the cardinality and capacity constraints.

This branch-and-prune approach is very fast, at most a few seconds, for very small $k$. So we embedded it into our branch-and-bound algorithm and run it at nodes where the remaining number of items to be filled in the current 
knapsack is very small (less or equal than $5$ in practice). To solve the original problem, we can also replace the global branch-and-bound method using this branch-and-prune approach for small initial values of $k$, in practice we choose $k \leq 10$.

\section{Numerical Results}\label{sec:numerical}
We coded the algorithm in C++. For the function evaluation (i.e., solving \eqref{sdp}) we implemented a
predictor-corrector variant of an interior point algorithm~\cite{HeReVaWo:96}. We use the ConicBundle
Library of Ch.~Helmberg~\cite{CB} as framework for the bundle method to solve~\eqref{sdp-met}.

We compare our method $(B\&C)$ to:
\begin{itemize}
 \item (Cplex): IBM CPLEX solver, version 12.6.2~\cite{CPL}, with default settings.
 \item (MIQCR+Cplex): our implementation of the MIQCR method~\cite{BEL11}. MIQCR uses a semidefinite relaxation in order to obtain a problem having a convexified objective function; the resulting convex integer problem can then be solved by standard solvers. We use the CSDP solver~\cite{Borchers2} for solving the semidefinite relaxation to convexify the objective function, and IBM CPLEX 12.6.2~\cite{CPL} with default settings 
       to solve the reformulated convex problem.
 \item (BiqCrunch): Also BiqCrunch~\cite{KrMaRo:14} is an algorithm based on semidefinite and polyhedral relaxations within a branch-and-bound framework. In BiqCrunch a quadratic regularization term is added to the objective function of the semidefinite problem and a quasi-Newton method is used to compute the bounds. We use the BiqCrunch solver enhanced with our primal and variable fixation heuristics described in Sect.~\ref{sec:heuristics}.
\end{itemize}

All experiments have been computed on an Intel i7-2600  quad core 3.4~GHz with 8~GB of RAM, using only one core.
The computational results have been obtained for randomly generated instances from~\cite{LET14} with up to 150~variables. The time limit for each approach is 3~hours.

In Table~\ref{tab:results} we display the run time of the overall algorithm, the gap at the root node, and the number of nodes produced during the branch-and-bound algorithm for each method. Each line of Table~\ref{tab:results} represents average values over 10~instances. We put the number of instances solved within the time limit into brackets in case not all 10~instances could be solved. Average values are computed only over instances solved within the time limit.

\begin{table}[h]																									
 \centering																									
 \tiny																									
 \setlength{\tabcolsep}{2.5pt}																									
\rotatebox{90}{
 \begin{tabular}{ccrrrrrrrrrrrrrrrrrrrrrrrrrr}																									
  \toprule																									
			&		&	&	\multicolumn{ 3}{c}{Cplex}	&	&	\multicolumn{ 3}{c}{MIQCR+Cplex}	&	&										
				\multicolumn{ 3}{c}{BiqCrunch}	&	&	\multicolumn{ 3}{c}{B\&C}	\\													
			&		&	&	Gap root \%	&	Time (s)	&	\#Nodes	&	&	Gap root \%	&	Time (s)	&	\#Nodes	&	&	
		Gap root \%	&	Time (s)	&	\#Nodes	&	&	Gap root \%	&	Time (s)	&	\#Nodes	\\
		n	&	$\delta$	&	&	&	&	&	&	&	&	&	&	&	&		&		&	&		&		\\
\cmidrule(lr){1-2} \cmidrule(lr){4-6}\cmidrule(lr){8-10}\cmidrule(lr){12-14} \cmidrule(lr){16-18}																									
{\em	50	}	&	{\em25}	&	&	102.7	&	3.7	&	3426.9	&	&	30.5	&	\textbf{1.0}	&	621.2	&	&	7.4	&	21.4	&	79.6	&	&	\textbf{0.9}	&	72.4	&	11.6	\\
			&	{\em50}	&	&	150.6	&	150.8	&	77807.9	&	&	25.2	&	\textbf{1.0}	&	1276.3	&	&	4.9	&	24.9	&	136.8	&	&	\textbf{1.3}	&	9.1	&	11.2	\\
			&	{\em75}	&	&	230.3	&	213.1	&	104419.5&	&	102.0	&	\textbf{0.7}	&	656.7	&	&	56.1	&	26.6	&	98.6	&	&	\textbf{0.6}	&	3.6	&	9.1	\\
			&	{\em100}&	&	356.5	&	(8) 53.1&	(8) 14228.8&	&	62.7	&	\textbf{1.5}	&	3620.0	&	&	31.4	&	23.0	&	89.6	&	&	\textbf{0.9}	&	73.0	&	38.1	\\
\cmidrule(lr){1-2} \cmidrule(lr){4-6}\cmidrule(lr){8-10}\cmidrule(lr){12-14} \cmidrule(lr){16-18}																									
{\em	60	}	&	{\em25}	&	&	60.8	&	3.0	&	917.1	&	&	127.4	&	\textbf{0.9}	&	621.2	&	&	123.0	&	32.2	&	85.2	&	&	\textbf{0.6}	&	18.3	&	18.4	\\
			&	{\em50}	&	&	93.7	&	282.4	&	134246.3&	&	15.1	&	\textbf{1.4}	&	1280.3	&	&	4.7	&	39.7	&	136.8	&	&	\textbf{2.0}	&	110.3	&	88.1	\\
			&	{\em75}	&	&	212.7	&	(9) 50.9&	(9) 8258.3&	&	137.5	&	\textbf{3.3}	&	7594.5	&	&	131.4	&	71.3	&	123.0	&	&	\textbf{1.3}	&	75.8	&	28.2	\\
			&	{\em100}&	&	284.5	&	(8) 188.6&	(8) 55411.8&	&	61.2	&	\textbf{3.2}	&	5808.1	&	&	47.8	&	63.0	&	147.2	&	&	\textbf{0.3}	&	21.5	&	18.4	\\
\cmidrule(lr){1-2} \cmidrule(lr){4-6}\cmidrule(lr){8-10}\cmidrule(lr){12-14} \cmidrule(lr){16-18}																									
{\em	70	}	&	{\em25}	&	&	130.2	&	23.7	&	12065.8	&	&	37.9	&	\textbf{3.4}	&	2884.9	&	&	13.8	&	109.6	&	147.7	&	&	\textbf{4.5}	&	259.3	&	42.0	\\
			&	{\em50}	&	&	177.1	&	(6) 213.8&	(6) 63859.7&	&	71.7	&	\textbf{8.4}	&	11221.8	&	&	59.2	&	141.0	&	207.4	&	&	\textbf{2.2}	&	128.2	&	139.7	\\
			&	{\em75}	&	&	382.4	&	(8) 873.0&	(8) 105465.6&	&	56.1	&	\textbf{16.2}	&	33821.0	&	&	17.4	&	196.2	&	211.7	&	&	\textbf{3.5}	&	246.8	&	114.9	\\
			&	{\em100}&	&	252.2	&	(4) 60.2&	(4) 10867.5&	&	59.6	&	\textbf{14.6}	&	25809.6	&	&	53.0	&	153.3	&	243.2	&	&	\textbf{4.0}	&	319.7	&	338.8	\\
\cmidrule(lr){1-2} \cmidrule(lr){4-6}\cmidrule(lr){8-10}\cmidrule(lr){12-14} \cmidrule(lr){16-18}																									
{\em	80	}	&	{\em25}	&	&	111.2	&	226.6	&	89013.4	&	&	33.5	&	\textbf{7.8}	&	6115.3	&	&	13.0	&	149.6	&	195.2	&	&	\textbf{7.5}	&	390.9	&	86.1	\\
			&	{\em50}	&	&	271.6	&	(8) 872.9&	(8) 181325.9&	&	55.0	&	\textbf{26.8}	&	36346.3	&	&	20.9	&	373.9	&	366.2	&	&	\textbf{8.6}	&	544.8	&	213.6	\\
			&	{\em75}	&	&	313.3	&	(5) 278.7&	(5) 14838.5&	&	82.0	&	\textbf{47.8}	&	96543.3	&	&	70.8	&	615.1	&	745.2	&	&	\textbf{2.6}	&	413.4	&	359.0	\\
			&	{\em100}&	&	473.0	&	(6) 1469.5&	(6) 98024.5&	&	43.0	&	\textbf{96.5}	&	216700.0&	&	17.1	&	717.5	&	804.6	&	&	\textbf{5.4}	&	1849.4	&	1219.7	\\
\cmidrule(lr){1-2} \cmidrule(lr){4-6}\cmidrule(lr){8-10}\cmidrule(lr){12-14} \cmidrule(lr){16-18}																									
{\em	90	}	&	{\em25}	&	&	118.5	&	(9) 585.9&	(9) 693035.0&	&	111.5	&	\textbf{23.3}	&	22836.9	&	&	107.3	&	188.6	&	390.4	&	&	\textbf{3.6}	&	430.6	&	94.1	\\
			&	{\em50}	&	&	248.6	&	(6) 3708.5&	(6) 312105.5&	&	82.2	&	\textbf{67.8}	&	99574.6	&	&	72.3	&	532.3	&	810.0	&	&	\textbf{3.4}	&	729.1	&	404.9	\\
			&	{\em75}	&	&	388.7	&	(2) 2850.5&	(2) 62190.5&	&	37.9	&	\textbf{735.1}	&	1348558.3&	&	14.2	&	1281.2	&	970.8	&	&	\textbf{8.7}	&	(7) 3234.1&	(7) 2233.1\\
			&	{\em100}&	&	390.0	&	(3) 146.2&	(3) 5047.5&	&	26.6	&	\textbf{180.4}	&	282966.1&	&	10.4	&	1094.7	&	5644.1	&	&	\textbf{6.5}	&	2740.9	&	1357.9	\\
\cmidrule(lr){1-2} \cmidrule(lr){4-6}\cmidrule(lr){8-10}\cmidrule(lr){12-14} \cmidrule(lr){16-18}																									
{\em	100	}	&	{\em25}	&	&	169.4	&	2308.1	&	623731.5&	&	74.4	&	\textbf{65.4}	&	71449.9	&	&	61.6	&	392.5	&	617.4	&	&	\textbf{10.9}	&	1583.1	&	284.0	\\
			&	{\em50}	&	&	145.7	&	(6) 1724.3&	(6) 122716.0&	&	17.5	&	\textbf{308.0}	&	465749.5&	&	\textbf{7.6}	&	986.9	&	882.0	&	&	8.0	&	(9) 3379.6&	(9) 1488.6\\
			&	{\em75}	&	&	270.9	&	(2) 4243.5&	(2) 88176.5&	&	21.8	&	\textbf{856.8}	&	1322350.0&	&	\textbf{6.8}	&	980.0	&	967.8	&	&	14.8	&	(7) 2613.0&	(7) 855.6\\
			&	{\em100}&	&	473.0	&	(5) 2658.8&	(5) 120959.0&	&	98.6	&	\textbf{649.7}	&	977246.9&	&	94.0	&	(9) 723.8&	(9) 5166.7&	&	\textbf{6.0}	&	(7) 318.1&	(7) 115.4\\
\cmidrule(lr){1-2} \cmidrule(lr){4-6}\cmidrule(lr){8-10}\cmidrule(lr){12-14} \cmidrule(lr){16-18}																									
{\em	110	}	&	{\em25}	&	&	124.0	&	(6) 277.4&	(6) 36270.2&	&	72.2	&	\textbf{327.0}	&	288129.4&	&	64.4	&	848.8	&	1003.6	&	&	\textbf{13.5}	&	(8) 2602.9&	(8) 810.4\\
			&	{\em50}	&	&	117.5	&	(3) 661.5&	(3) 55327.0&	&	14.3	&	1188.1	&	1089556.0&	&	\textbf{4.7}	&	\textbf{1010.2}	&	727.8	&	&	5.7	&	(7) 2065.1&	(7) 652.3\\
			&	{\em75}	&	&	580.7	&	(6) 908.8	&	(6) 35891.8&	&	138.7	&	(7) 27.4&	(7) 37408.8&	&	118.8	&	\textbf{(8) 2305.7}&	(8) 4523.3&	&	\textbf{10.7}	&	(7) 1062.7&	(7) 297.7\\
			&	{\em100} &	&	332.2	&	(1) 1911.6&	(1) 118552.0&	&	19.6	&	\textbf{(8) 758.0}&	(8) 956886.3&	&	7.1	&	(8) 1438.1&	(8) 1575.0&	&	\textbf{7.0}	&	(6) 1789.2&	(6) 511.7\\
\cmidrule(lr){1-2} \cmidrule(lr){4-6}\cmidrule(lr){8-10}\cmidrule(lr){12-14} \cmidrule(lr){16-18}																									
{\em	120	}	&	{\em25}	&	&	55.1	&	(6) 320.1&	(6) 94936.5&	&	95.3	&	1771.4	&	1644176.9&	&	111.8	&	\textbf{424.3}	&	447.6	&	&	\textbf{5.7}	&	(8) 1872.3&	(8) 317.3\\
			&	{\em50}	&	&	288.1	&	(3) 2995.9&	(3) 81429.7&	&	90.3	&	(7) 1888.9&	(7) 1554792.4&	&	82.5	&	\textbf{(8) 1073.9}&	(8) 821	&	&	\textbf{10.2}	&	(6) 1725.3&	(6) 427.5\\
			&	{\em75}	&	&	507.6	&	(5) 305.0	&	(5) 11101.2&	&	133.8	&	(6) 484.9&	(6) 177043.8&	&	128.7	&	\textbf{(9) 3001.6}&	(9) 7996.3&	&	\textbf{7.9}	&	(5) 0.8	&	(5) 0.0	\\
			&	{\em100} &	&	179.7	&	(3) 41.9	&	(3) 4166.5	&	&	66.2	&	(6) 61.6&	(6) 36075.0&	&	68.7	&	\textbf{(9) 1552.6}&	(9) 969.0&	&	\textbf{3.6}	&	(6) 683.4&	(6) 144.4\\
\cmidrule(lr){1-2} \cmidrule(lr){4-6}\cmidrule(lr){8-10}\cmidrule(lr){12-14} \cmidrule(lr){16-18}																									
{\em	130	}	&	{\em25}	&	&	129.1	&	(6) 2014.5&	(6) 383586.2&	&	24.7	&	\textbf{1256.6}	&	698989.0	&		&	10.6	&	3341.8	&	2520.8	&	&	\textbf{7.1}	&	(8) 4194.0&	(8) 850.7\\
			&	{\em50}	&	&	411.8	&	(4) 3246.9	&	(4) 245787.0&	&	67.3	&	(7) 493.0	&	(7) 384516.5	&	&	50.5	&	\textbf{(8) 1719.6}	&	(8) 2330.7&	&	\textbf{6.9}	&	(7) 2590.8&	(7) 450.8\\
			&	{\em75}	&	&	207.3	&	(0)	&	(0)	&	&	12.2	&	(5) 4975.3	&	(5) 3138617.2&	&	\textbf{3.8}	&	\textbf{(9) 2630.0}	&	(9) 1000.3	&	&	11.9	&	(2) 6813.5&	(2) 1430.0\\
			&	{\em100} &	&	383.5	&	(0)	&	(0)	&	&	21.1	&	(4) 2250.4	&	(4) 1170285.3&	&	\textbf{8.8}	&	\textbf{(6) 4365.0}	&	(6) 2437.4	&	&	14.1	&	(3) 2012.0&	(3) 83.5\\
\cmidrule(lr){1-2} \cmidrule(lr){4-6}\cmidrule(lr){8-10}\cmidrule(lr){12-14} \cmidrule(lr){16-18}																									
{\em	140	}	&	{\em25}	&	&	207.9	&	(5) 15.0&	(5) 1180.5&	&	48.8	&	(8) 1770.8	&	(8) 654605.8&	&	44.4	&	\textbf{2624.2}	&	1993.3	&	&	\textbf{13.3}	&	(4) 2561.8	&	(4) 298.5	\\
			&	{\em50}	&	&	306.2	&	(1) 2401.8& (1) 106348.0&	&	36.8	&	(5) 2692.3	&	(5) 2306370.0&	&	\textbf{16.4}	&	\textbf{(7) 4809.5}	&	(7) 1134.1	&	&	24.0	&	(3) 3360.3&	(3) 213.3\\
			&	{\em75}	&	&	259.0	&	(0)	&	(0)	&	&	19.6	&	(4) 2263.5	&	(4) 1520163.5&	&	\textbf{8.0}	&	\textbf{(5) 4065.2}	&	(5) 1773.4	&	&	12.5	&	(1) 431.0&	(1) 0.0	\\
			&	{\em100} &	&	647.8	&	(2) 64.1	&	(2) 483.0	&	&	49.4	&	(4) 1042.9	&	(4) 1238929.3&	&	37.1	&	\textbf{(6) 6123.7}	&	(6) 17745.3	&	&	\textbf{13.3}	&	(4) 2561.7&	(4) 298.5\\
\cmidrule(lr){1-2} \cmidrule(lr){4-6}\cmidrule(lr){8-10}\cmidrule(lr){12-14} \cmidrule(lr){16-18}																									
{\em	150	}	&	{\em25}	&	&	103.7	&	(4) 1744.3&	(4) 202004.0&	&	67.4	&	(6) 587.1	&	(6) 552495.8&	&	69.1	&	\textbf{(8) 3203.9}	&	(8) 994.8&	&	\textbf{12.6}	&	(3) 1458.0	 &	(3) 164.3	\\
			&	{\em50}	&	&	105.6	&	(5) 91.8	&	(5) 5591.7	&	&	98.2	&	\textbf{(7) 2240.0}	&	(7) 935027.4&	&	101.4	&	(7) 2761.5&	(7) 2349.4&	&	\textbf{8.8}	&	(5) 2797.7	&	(5) 3453.7	\\
			&	{\em75}	&	&	496.9	&	(0)	&	(0)	&	&	\textbf{7.9}	&	\textbf{(5) 957.2}	&	(5) 300455.4	&	&	23.1	&	(3) 3908.7	&	(3) 876.3	&	&	17.5	&	(2) 155.0	&	(2) 0.0	\\
			&	{\em100} &	&	171.1	&	(1) 1039.3	&	(1) 24546.0	&	&	43.7	&	(3) 3493.0	&	(3) 5264462.0	&	&	\textbf{3.3}	&	\textbf{(4) 4320.1}	&	(4) 2171.0	&	&	8.1	&	(3) 5391.7	&	(3) 554.3	\\
\cmidrule(lr){1-2} \cmidrule(lr){4-6}\cmidrule(lr){8-10}\cmidrule(lr){12-14} \cmidrule(lr){16-18}																									
		\multicolumn{ 2}{c}{{\em	Avg	}}			&	&	258.5	&	(236) 787.0&	(236)127524.4&	&	59.0	&	(372) 570.0&	(372) 902502.7&	&	45.7	&	\textbf{(394) 1215.4}	&(394)	1487.2	&	&	\textbf{7.4}	&	(328) 1250.5&	(328) 433.8\\
    \bottomrule																									
\end{tabular}}																									
 \caption{Numerical results comparing four approaches. (Time limit: 3~hours)} \label{tab:results}																									
\end{table}

The numerical experiments demonstrate that the methods having semidefinite optimization inside clearly outperform Cplex. In fact, Cplex already fails to solve all instances of size $n=50$ within the time limit.

Instances with up to $n=100$ variables can be solved most efficiently by the MIQCR approach, i.e., finding a convexified problem via semidefinite optimization and then solve the resulting convex problem using Cplex. 

For $n>100$, BiqCrunch performs best in terms of overall run time, but the domincance to MIQCR and our approach is not significant. 

Our new approach provides by far the smallest gap at the root node. The high quality of our bound is also reflected in the number of nodes in the branch-and-bound tree. Our method explores a substantial smaller number of nodes than the other approaches.

Our approach is not superior to MIQCR or BiqCrunch in terms of overall computation time, however, the implementation is a prototype and there is room for speeding up the approach by experimenting with different settings in the branch-and-bound framework (such as branching strategies) as well as parameter settings in the bundle algorithm and in the update of the set of triangle inequalities. This is currently under investigation. 

\section{Conclusion}\label{sec:conclusion}
The 0-1 $k$-item quadratic knapsack problem is a challenging problem, as it includes two NP-hard problems, namely quadratic knapsack and $k$-cluster. We review approaches to solve this problem to optimality and introduce a new method, where the bound computation is based on a semidefinite relaxation. The derived basic semidefinite relaxation has only simple constraints, in fact all constraints are of rank one. This can be exploited in interior point methods to efficiently compute the system matrix. We strengthen the relaxation using triangle inequalities and solve the resulting semidefinite problem by a dynamic version of the bundle method. 

To have a comparison with state of the art algorithms we implement the convexification algorithm MIQCR~\cite{BEL11}, use BiqCrunch~\cite{KrMaRo:14} enhanced with our primal heuristics, and run Cplex. The numerical results prove that CPLEX is clearly outperformed by all the methods based on semidefinite programming. Our new method provides the tightest bound at the root node, while the overall computation time is smallest for MIQCR for $n \le 100$ and BiqCrunch for larger $n$. An optimized implementation and a study of the best parameter settings for the various components inside our code is subject of further study.

\bibliographystyle{plain}
\bibliography{../biblio}

\end{document}